\documentclass[a4paper]{article}
\usepackage{amssymb}
\newcommand{\bea}{\begin{eqnarray}}
\newcommand{\eea}{\end{eqnarray}}
\newcommand{\ba}{\begin{array}}
\newcommand{\ea}{\end{array}}
\newcommand{\bc}{\begin{center}}
\newcommand{\ec}{\end{center}}

\newcommand{\be}{\begin{equation}}
\newcommand{\edc}{\end{document}}
\newcommand{\ee}{\end{equation}}

\newcommand{\dsf}{\displaystyle\frac}

\def\m{\mu}

\def\g{\gamma}

\def\l{\lambda}
\def\g{\gamma}
\def\r{\rho}

\def\Q{\mathbb{Q}}
\def\Z{\mathbb{Z}}
\def\R{\mathbb{R}}
\def\N{\mathbb{N}}
\def\C{\mathbb{C}}
\begin{document}
\bc {\Large {\bf ON A CLASS OF RATIONAL $P$-ADIC DYNAMICAL SYSTEMS}}\\[6mm]

{\large MUROD KHAMRAEV\\
{\it Department of Mathematics\\
University of Manchester\\
Manchester, UK\\
e-mail: {\tt Murod.Khamraev@postgrad.manchester.ac.uk}}\\[1cm]
FARRUKH MUKHAMEDOV\\
{\it
Department of Mechanics and Mathematics\\
National University of Uzbekistan\\
Vuzgorodok, Tashkent, 700174, Uzbekistan,\\
e-amil: {\tt far75m@yandex.ru}}}\\[3mm]

\ec %\Large
\begin{abstract}
In this paper we investigate the behavior of trajectories of one
class of rational $p$-adic dynamical systems in complex $p$-adic
field $\C_p$. We studied Siegel disks and attractors of such
dynamical systems. We found the basin of the attractor of the
system. It is proved that such dynamical systems are not ergodic
on a unit sphere with respect to the Haar measure.
\\[2mm]
%\ec
\end{abstract}
%\large

\section{Introduction}

The $p$-adic numbers were first introduced by the German
mathematician K.Hensel. During a century after the discovery of
$p$-adic numbers, they were considered mainly objects of pure
mathematics. Starting from 1980's various models described in the
language of $p$-adic analysis have been actively studied. More
precisely, some models over the field of $p$-adic numbers have
been considered, which is due to the assumption that $p$-adic
numbers provide a more exact and more adequate description of
micro-world phenomena. Numerous applications of these numbers to
theoretical physics have been proposed in papers \cite{ADFV, FW,
MP, V1, V2} to quantum mechanics \cite{Kh1}, to $p$-adic - valued
physical observable \cite{Kh1} and many others \cite{Kh2, VVZ}.

The study of $p$-adic dynamical systems arises in Diophantine
geometry in the constructions of canonical heights, used for
counting rational points on algebraic vertices over a number
field, as in \cite{CS}. In \cite{Kh4, TVW} $p$-adic field have
arisen in physics in the theory of superstrings, promoting
questions about their dynamics. Also some applications of $p$-adic
dynamical systems to some biological, physical systems have been
proposed in \cite{ABKO, AKK1, AKK2, DGKS, Kh5, KhN}. Other studies
of non-Archimedean dynamics in the neighborhood of a periodic and
of the counting of periodic points over global fields using local
fields appear in \cite{HY, L, P}. It is known, that analytic
functions play important role in complex analysis. In the $p$-adic
analysis the rational functions play a similar role to the
analytic functions in complex analysis \cite{E, R}. Therefore,
naturally there arises a need to study the dynamics of these
functions in the $p$-adic analysis. On the other hand, such
$p$-adic dynamical systems appear in the process of studying
$p$-adic Gibbs measures \cite{gmr, KM, KMR, MR1}. In \cite{B1, B2}
dynamics on the Fatou set of a rational function defined over some
finite extension of $\Q_p$ has been studied, besides, an analogue
of Sullivan's no wandering domains' theorem for $p$-adic rational
Functions, which have no wild recurrent Julia critical points,
were proved. Recently, in \cite{MR2} some linear-rational
dynamical systems on complex $p$-adic numbers $\C_p$ have been
investigated. In \cite{AKTS} the behavior  of a $p$-adic dynamical
system $f(x)=x^n$ in the fields of $p$-adic numbers $\Q_p$ and
$\C_p$ was investigated. Some ergodic properties of that dynamical
system have been considered in \cite{GKL}. In \cite{TVW, AKK2}
certain limit behaviors of a $p$-adic dynamical system of the form
$f(x)=x^2+c$, $c\in\Q_p$ were investigated on $\Q_p$. In
\cite{GKL} it was shown that such a dynamical system is ergodic at
$c=0$. In \cite{Sh} the Fatou and Julia sets were found of a
$p$-adic dynamical system $f(x)=\l x(x-1)$, $\l\in\Q_p$. Using
them a topological entropy one is also calculated. These
investigations lead us to consideration of a rational perturbation
of the dynamical system $f(x)=x^2$. In the paper we will
investigate the behavior of trajectory of a rational $p$-adic
dynamical systems of the form $f(x)=\frac{ax^{2}}{bx+1}$ in
$\C_p$. We will study Siegel disks and attractors of such
dynamical systems. Besides, it is shown that the considered
dynamical system is not ergodic with respect to the Haar measure
on a sphere. Note the basics of $p$-adic analysis, $p$-adic
mathematical physics are explained in \cite{G, Ko, VVZ}.

\section{Preliminaries}

\subsection{$p$-adic numbers}

Let $\Q$ be the field of rational numbers. The greatest common
divisor of the positive integers $n$ and $m$ is denoted by
$(n,m)$. Every rational number $x\neq 0$ can be represented in the
form $x=p^r\dsf{n}{m}$, where $r,n\in\mathbb{Z}$, $m$ is a
positive integer, $(p,n)=1$, $(p,m)=1$ and $p$ is a fixed prime
number. The $p$-adic norm of $x$ is given by
$$
|x|_p=\left\{
\ba{ll}
p^{-r} & \ \textrm{ for $x\neq 0$}\\
0 &\ \textrm{ for $x=0$}.\\
\ea
\right.
$$

It has the following properties:

1) $|x|_p\geq 0$ and $|x|_p=0$ if and only if $x=0$,

2) $|xy|_p=|x|_p|y|_p$,

3) the strong triangle inequality
$$
|x+y|_p\leq\max\{|x|_p,|y|_p\},
$$

3.1) if $|x|_p\neq |y|_p$ then $|x-y|_p=\max\{|x|_p,|y|_p\}$,

3.2) if $|x|_p=|y|_p$ then $|x-y|_p\leq |x|_p$,

this is a non-Archimedean inequality.

The completion of $\Q$ with respect to $p$-adic norm defines
the $p$-adic field which is denoted by $\Q_p$.

Any $p$-adic number $x\neq 0$ can be uniquely represented in the
canonical series:
$$
x=p^{\g(x)}(x_0+x_1p+x_2p^2+...) ,
\eqno(2.1)
$$
where $\g=\g(x)\in\Z$ and $x_j$ are integers, $0\leq x_j\leq p-1$, $x_0>0$,
$j=0,1,2,...$ (see more detail \cite{G},\cite{Ko}).
Observe that in this case $|x|_p=p^{-\g(x)}$.

The algebraic completion of $\Q_p$ is denoted by $\C_p$ and it is called
{\it complex $p$-adic numbers}.  For any $a\in\C_p$ and $r>0$ denote
$$
\bar B_r(a)=\{x\in\C_p : |x-a|_p\leq r\},\ \ B_r(a)=\{x\in\C_p : |x-a|_p< r\},
$$
$$
S_r(a)=\{x\in\C_p : |x-a|_p= r\}.
$$

A function $f:B_r(a)\to\C_p$ is said to be {\it analytic} if it
can be represented by
$$
f(x)=\sum_{n=0}^{\infty}f_n(x-a)^n, \ \ \ f_n\in \C_p,
$$ which converges uniformly on the ball $B_r(a)$.
For a review of basic properties of analytic functions we refer to
\cite{E}.

\subsection{Dynamical systems in $\C_p$}

In this section we recall some known facts about dynamical systems
$(f,B)$ in $\C_p$, where $f: x\in B\to f(x)\in B$ is an analytic
function and $B=B_r(a)$ or $\C_p$.

Recall some standard terminology of the theory of dynamical
systems (see for example \cite{PJS}). Let $f:B\to B$ be an
analytic function. Denote $x^{(n)}=f^n(x^{(0)})$, where $x^0\in B$
and $f^n(x)=\underbrace{f\circ\dots\circ f(x)}_n$.
 If $f(x^{(0)})=x^{(0)}$ then $x^{(0)}$
is called a {\it fixed point}. A fixed point $x^{(0)}$ is called
an {\it attractor} if there exists a neighborhood $U(x^{(0)})$ of
$x^{(0)}$ such that for all points $y\in U(x^{(0)})$ it holds
$\lim\limits_{n\to\infty}y^{(n)}=x^{(0)}$, where $y^{(n)}=f^n(y)$.
If $x^{(0)}$ is an attractor then its {\it basin of attraction} is
$$
A(x^{(0)})=\{y\in \C_p :\ y^{(n)}\to x^{(0)}, \ n\to\infty\}.
$$
A fixed point $x^{(0)}$ is called {\it repeller} if there  exists
a neighborhood $U(x^{(0)})$ of $x^{(0)}$ such that
$|f(x)-x^{(0)}|_p>|x-x^{(0)}|_p$ for $x\in U(x^{(0)})$, $x\neq
x^{(0)}$. Let $x^{(0)}$ be a fixed point of a function $f(x)$. The
ball $B_r(x^{(0)})$ (contained in $B$) is said to be a {\it Siegel
disk} if each sphere $S_{\r}(x^{(0)})$, $\r<r$ is an invariant
sphere of $f(x)$, i.e. if $x\in S_{\r}(x^{(0)})$ then all iterated
points $x^{(n)}\in S_{\r}(x^{(0)})$ for all $n=1,2\dots$. The
union of all Siegel disks with the center at $x^{(0)}$ is said to
{\it a maximum Siegel disk} and is denoted by $SI(x^{(0)})$.

{\bf Remark 2.1.}\cite{AKTS} In complex geometry, the center of a
disk is uniquely determined by the disk, and different fixed
points cannot have the same Siegel disks. In non-Archimedean
geometry, a center of a disk is nothing but a point which belongs
to the disk. Therefore, in principle, different fixed points may
have the same Siegel disk.

Let $x^{(0)}$ be a fixed point of an analytic function $f(x)$. Set
$$
\l=\frac{d}{dx}f(x^{(0)}).
$$

The point $x^{(0)}$ is called {\it attractive} if $0\leq |\l|_p<1$, {\it indifferent} if
$|\l|_p=1$, and {\it repelling} if $|\l|_p>1$.

{\bf Theorem 2.1.}\cite{AKTS} {\it Let $x^{(0)}$ be a fixed point
of an analytic function $f:B\to B$. The following assertions hold

1. if $x^{(0)}$ is an attractive point of $f$, then it is an
attractor of the dynamical system $(f,B)$. If $r>0$ satisfies the
inequality
$$
q=\max_{1\leq n<\infty}\bigg|\frac{1}{n!}\frac{d^nf}{dx^n}(x^{(0)})\bigg|_pr^{n-1}<1
\eqno(2.2)
$$
and $B_r(x^{(0)})\subset B$ then $B_r(x^{(0)})\subset A(x^{(0)})$;

2. if $x^{(0)}$ is an indifferent point of $f$ then it is the
center of a Siegel disk. If $r$ satisfies the inequality
$$
s=\max_{2\leq
n<\infty}\bigg|\frac{1}{n!}\frac{d^nf}{dx^n}(x^{(0)})\bigg|_pr^{n-1}<1
\eqno(2.3)
$$
and $B_r(x^{(0)})\subset B$ then $B_r(x^{(0)})\subset
SI(x^{(0)})$;

3. if $x^{(0)}$ is a repelling point of $f$ then $x^{(0)}$ is a
repeller of the dynamical system $(f,B)$.}

\section{$p$-adic dynamical system of the form
$f(x)=\frac{ax^{2}}{bx+1}$}

In this section we consider dynamical system associated with the function
$f:\C_p\to\C_p$ defined by
$$
f(x)=\frac{ax^{2}}{bx+1}, \ \ \ \ a,b\in\mathbb{C}_{p},\ a,b\neq 0,\ a\neq b
\eqno(3.1)
$$
where $x\neq -\dsf{1}{b}$. Denote $P=-\dsf{1}{b}$ and $D=\C_p\setminus\{P\}$.

It is not difficult to check that fixed points of the function (3.1) are
$$
x_{1}=0\mbox{
and } \ \ x_{2}=\frac{1}{a-b}.
\eqno(3.2)
$$
Define the following sets
$$
\Omega=\left\{x\in D:\exists i\geq
1\Rightarrow x^{(i)}=P \right\},
$$
$$
\Psi=\left\{x\in
D\backslash\{x_{2}\}:\exists j\geq 1\Rightarrow
x^{(j)}=x_{2}\right\},
$$
$$
\Sigma=\left\{x\in D\backslash\{x_{1},x_{2}\}:\exists
k\geq 2 \Rightarrow x^{(k)}=x\right\}.
$$
From (3.1) we can easily prove the following assertions.

{\bf Proposition 3.1.} {\it The following is true:}
\begin{enumerate}
  \item \textit{$\Omega\neq\emptyset,$ $\Psi\neq\emptyset,$ and
  $\Sigma\neq\emptyset.$}
  \item \textit{$\Omega\bigcap\Psi=\emptyset,$ $\Omega\bigcap\Sigma=\emptyset,$
  and $\Psi\bigcap\Sigma=\emptyset.$}
\end{enumerate}

{\bf Proposition 3.2.} \textit{The following equalities hold for
$\forall x \in D$:}
\begin{enumerate}
  \item[$\Delta 1$] $f(x)(x-P)=\left(\dsf{a}{b}\right)x^{2}.$
  \item[$\Delta 2$]
  $(f(x)-x)(x-P)=\left(\dsf{a-b}{b}\right)x(x-x_{2}).$
  \item[$\Delta 3$]
  $(f(x)-x_{2})(x-P)=\dsf{a}{b}\left(x+\dsf{1}{a}\right)(x-x_{2}).$
\end{enumerate}

{\bf Remark 3.1.} Sometimes it is more convenient to use
the last equality in a slightly different form - simple algebraic
operations allow us to do that - namely:
$$(f(x)-x_{2})\left(1+\dsf{(a-b)b}{a}(x-x_{2})\right)
=(a-b)(x-x_{2})\left(x+\frac{1}{a}\right).$$

Let us calculate the value of $\frac{d^{n}}{d^{n}x}f(x)$ for
$n\geq 1$ in $x_{1}=0$ and $x_{2}=\frac{1}{a-b}$. Note, that
$f(x)=\frac{a}{b}\left(x-\frac{1}{b}+\frac{1}{b^{2}}\left(\frac{1}{x+\frac{1}{b}}\right)\right)$,
thus
$$f'(x)=\frac{a}{b}\left(1-\frac{1}{b^{2}}\frac{1}{\left(x+\frac{1}{b}\right)^{2}}\right),$$
and for $n\geq 2$ we have
$$\frac{d^{n}}{d^{n}x}f(x)=\frac{a}{b^{3}}(-1)^{n}n!\frac{1}{\left(x+\frac{1}{b}\right)^{n+1}}.$$
Hence, we have that $f'(x_{1})=0$ and $f'(x_{2})=\frac{2a-b}{a}$
and for $n\geq 2$ we deduce that
$$\frac{d^{n}}{d^{n}x}f(x_{1})=(-1)^{n}n!ab^{n-2}$$ and $$\frac{d^{n}}{d^{n}x}f(x_{2})=(-1)^{n}n!a^{-n}b^{n-2}(a-b)^{n+1}.$$
According to Theorem 2.1 we have, that $x_{1}=0$ is the
attracting point of the dynamical system $(f(x),D)$,
since $|f'(x_{1})|_p=0<1$. For $x_{2}=\frac{1}{a-b}$ its nature is
defined by the magnitude of the value $\frac{|2a-b|_p}{|a|_p}$.

Analyzing all possible values of $|2a|_p$ and $|b|_p$ three main possibilities arise:

\begin{enumerate}
\item $x_{2}$ is the repelling point $\Leftrightarrow$ $|2a-b|_p>|a|_p$.
This condition splits to the following one:
\begin{enumerate}
\item $|a|_p<|b|_p$ for arbitrary prime $p\geq 2$.
\end{enumerate}

\item $x_{2}$ is the indifferent point $\Leftrightarrow$ $|2a-b|_p=|a|_p$. This condition
splits to the following ones:
\begin{enumerate}
 \item $|2a|_p>|b|_p$  and $\forall p>2$;

 \item $|2a|_p=|b|_p$, $\forall p>2$ and $|2a-b|_p=|a|_p$;

 \item $|a|_p=|b|_p$ and $p=2$.
\end{enumerate}

\item $x_{2}$ is the attracting point $\Leftrightarrow$ $|2a-b|_p<|a|_p$. This condition
splits to the following ones:

    \begin{enumerate}
      \item $|2a|_p>|b|_p$ and $p=2$;

      \item $|2a|_p=|b|_p$  and $\|2a-b\|<|a|_p$ for arbitrary $p\geq 3$ or $p=2$ and
      $|2a|_p=|b|_p$;

      \item  $|2a|_p<|b|_p<|a|_p$ and $p=2$.
    \end{enumerate}
\end{enumerate}

We consider these possibilities in succession. The main points of
interest are the size of the attractor $A(x_{1})$ of the
attracting point $x_{1}=0$, and the sizes of $A(x_{2})$ and
$SI(x_{2})$ when the point $x_{2}=\frac{1}{a-b}$ is the attracting
and the indifferent point respectively.

To estimate the size of the $A(x_{1})$ we make use of the Theorem
2.1 - the requirement $$\max_{1\leq
n<\infty}\left|\frac{1}{n!}\frac{d}{dx}f(x_1)\right|_pr^{n-1}<1$$
takes the form of the following inequality in this case:
$$\max_{n\geq 1}\frac{|a|_p}{|b|_p}(|b|_pr)^{n-1}<1.$$ Multiplying
both sides of this inequality by $\frac{|b|_p}{|a|_p}$ and
recalling that $|f(x_{1})|_p=0$ we obtain the inequality
$$\max_{n\geq 1}(|b|_pr)^{n-1}=\max_{n\geq 2}(|b|_pr)^{n-1}<\frac{|b|_p}{|a|_p}.$$

Denote this condition by $\Gamma 1.$

The same theorem helps us find the radius of the Siegel disk of
$x_{2}$ when $x_{2}$ is an indifferent point - the requirement
$$\max_{2\leq
n<\infty}\left|\frac{1}{n!}\frac{d}{dx}f(x_2)\right|_pr^{n-1}<1$$
takes the form $$\max_{n\geq
2}\frac{|b|_p^{n-2}}{|a|_p^{n}}r^{n-1}|a-b|_p^{n+1}<1,$$ that is
more convenient to use in this case. Denote this condition by
$\Gamma 2.$

When $x_{2}$ is the attracting point, the size of the attractor
$A(x_{2})$ can be estimated using the inequality required by the
Theorem 2.1 - it takes the following form then -
$$\max_{n\geq
1}\frac{|b|_p^{n-2}}{|a|_p^{n}}r^{n-1}|a-b|_p^{n+1}<1.$$ This
condition we denote by $\Gamma 3.$

\subsection{$x_2$ is repelling point}

As was demonstrated, this is only possible in the case of
$|a|_p<|b|_p$. Everywhere in this subsection we will assume that
$|a|_p<|b|_p$.

Define for $\forall n\geq 0$,
$$
r_{n}=\frac{|b|_p^{n-1}}{|a|_p^{n}}, \ \ \ l_{n+1}=\frac{|a|_p^{n}}{|b|_p^{n+1}} \ \ \ \ l_0=0.
$$
Note, that $r_{0}=\dsf{1}{|b|_p}$ and
$r_{1}=\dsf{1}{|a|_p}$ - these values will be used extensively in
what follows.

It is easy to compute that the  $\Gamma 1$ can be satisfied only if $r\leq
r_{0}$ and then, if $r=r_{0}$ the maximal value is attained for
$\forall n\in\mathbb{N}$ and is equal to $\dsf{|a|_p}{|b|_p}<1$,
whereas, for $r<r_{0}$ the maximal value is attained for $n=2$ and
is equal to $\dsf{|a|_p}{|b|_p}<1$. Therefore, since
$B_{r_{0}}(x_{1})\subseteq D$, we  now know that
$B_{r_{0}}(x_{1})\subseteq A(x_{1})$. It is interesting to note,
that necessarily $x_{2}\notin A(x_{1})$, but
$|x_{2}|_p=r_{0}$,  and
thus, $x_{2}\in S_{r_{0}}(x_{1})\subset\bar{B}_{r_{0}}(x_{1})$.
So, $S_{r_{0}}(x_{1})\nsubseteq A(x_{1})$.

Now, suppose that $x\in D$ and $|x|_p>r_{0}$ then, since
$|P|_p=r_{0}$ the strong triangle inequality and the equality
3.1 with $\Delta 1$ imply that
$$
|f(x)|_p=\dsf{|a|_p}{|b|_p}|x|_p.
\eqno(3.3)
$$
Let us assume that $|x|\neq r_{n}$ for $\forall
n\geq 1$ and $|x|_p>r_{0}$, then $\exists j\in\mathbb{N}$, such
that $r_{j-1}<|x|_p<r_{j}$. Hence, using (3.3) we have
$r_{j-2}<|x^{(1)}|_p=\dsf{|a|_p}{|b|_p}|x|_p<r_{j-1}$, and thus, if
we continue iterating, we get the following estimate
$|x^{(j-1)}|_p<r_{0}$. i.e. $x^{(j-1)}\in B_{r_{0}}(x_{1})$. This
means, that for $\forall x\in D$ with $x\notin S_{r_{n}}(x_1)$ for
$\forall n\geq 0$, we find $x\in A(x_{1})$ or
$D\backslash\left(\bigcup\limits_{i=1}^{\infty}S_{r_{i}}(x_{1})\right)\subseteq
A(x_{1})$.

It remains to consider a case $|x|_p=r_0$. It is clear that
$|x-P|_p\leq r_0$, so assume that $|x-P|_p\neq l_n$ for all $n\geq 1$,
then $\exists j\in\mathbb{N}$, such
that $l_{j+1}<|x-P|_p<l_{j}$. Whence from $\Delta 1$ we infer
$$
r_{j-2}<|f(x)|_p<r_{j-1}.
$$
Hence, as said above, we obtain that $f(x)\in A(x_1)$ this means
$x\in A(x_1)$. Consequently, we have
$$
D\backslash\left(\bigcup_{i=1}^{\infty}S_{r_{i}}(x_{1})\cup
\bigcup_{j=0}^{\infty}S_{l_{j}}(P)\right)\subseteq
A(x_{1}).
$$

Thus, we have proved the following

{\bf Theorem 3.3} \textit{For the dynamical system
$\left(\frac{ax^{2}}{bx+1},D\right)$ the set of attraction of the
attracting fixed point $x_{1}$ is such, that
$$
D\backslash\left(\bigcup_{i=1}^{\infty}S_{r_{i}}(x_{1})\cup
\bigcup_{j=0}^{\infty}S_{l_{j}}(P)\right)\subseteq
A(x_{1}).
$$}

\textbf{Corollary 3.4.} \textit{For the dynamical system
$\left(\frac{ax^{2}}{bx+1},D\right)$ sets $\Omega,\Psi,\Sigma$ are
such, that
$$\Omega\bigcup\Psi\bigcup\Sigma\subset\left(\bigcup_{i=1}^{\infty}S_{r_{i}}(x_{1})\cup
\bigcup_{j=0}^{\infty}S_{l_{j}}(P)\right).$$}

Denote
\bea
B=\{x\in D : \ \forall n\in\N \ \ \exists r_n \ \textrm{or} \ l_n \ \textrm{such that} \nonumber \\
|f^{(n)}(x)|_p=r_n \ \textrm{or} \ \ |f^{(n)}(x)-P|_p=l_n\}.\nonumber
\eea

{\bf Corollary 3.5} \textit{For the dynamical system
$\left(\frac{ax^{2}}{bx+1},D\right)$ the set of attraction of the
attracting fixed point $x_{1}$ is such, that
$$
D\backslash B=A(x_{1}).
$$}

\subsection{$x_2$ is an indifferent point}

Here, as was demonstrated, we have three distinct cases, namely:

\begin{enumerate}
    \item The case of arbitrary $p>2$ and $|2a|_p>|b|_p$.
    \item The case of arbitrary $p>2$, $|2a|_p=|b|_p$ and $|2a-b|_p=|a|_p$.
    \item The case of $p=2$ and $|b|_p=|a|_p$.
\end{enumerate}

Now consider these cases separately.

{\bf Case of $p>2$ and $|2a|_p>|b|_p$}.

In this case $|2a|_p=|a|_p>|b|_p$,
$\dsf{1}{|a-b|_p}=\dsf{1}{|a|_p}=|x_{2}|_p$, and
$\dsf{1}{|a|_p}<\dsf{1}{|b|_p}=|P|_p$. The verification of the
condition $\Gamma 1$ shows that for all $r>0$ such that $r<r_1$ we
have $B_{r}(x_1)\subseteq A(x_{1})$.

Let us suppose that $|x|_p\geq r_1$, then
$\Delta 1$ with $|x-P|_p\leq\max\bigg\{|x|_p,\dsf{1}{|b|_p}\bigg\}$ implies that
$$
|f(x)|_p\geq\bigg|\dsf{a}{b}\bigg|_p\dsf{|x|_p^2}{\max\{|x|_p,\frac{1}{|b|_p}\}}\geq |x|_p
$$
which means $x$ does not lie in the $A(x_{1})$. Hence
 $B_{r_{1}}(x_{1})=A(x_{1})$.

Now, consider $x_{2}$ and its Siegel disk. The condition for $r$
is $\Gamma 2$ and it suggests with $|2a|_p>|b|_p$ that
$$
\max_{n\geq 2}\frac{|a|_p}{|b|_p}\left(|b|_pr\right)^{n-1}<1.
$$
From this we easily find that for every $0<r<r_{1}$ we have
$ B_{r}(x_{2})\subset SI(x_{2})$. And, therefore,
$B_{r_{1}}(x_{2})\subseteq SI(x_{2})$.

It is clear that  $-\frac{1}{a}=x_{0}\in \bar{B}_{r_{1}}(x_{2})$,
since
$$
\left|-\frac{1}{a}-x_{2}\right|_p=\left|\frac{2a-b}{a(a-b)}\right|_p=
\frac{1}{|a|_p}=r_{1},
$$
but evidently $-\dsf{1}{a}\notin SI(x_{2})$, since
$f(x_0)=x_2$. Hence we have precisely
$B_{r_{1}}(x_{2})=SI(x_{2})$.

This proves the following

\textbf{Theorem 3.6.} \textit{When $\forall p>2$ and
$|2a|_p>|b|_p$ the points $x_{1}$ and $x_{2}$ are respectively the
attracting point and the indifferent point of the dynamical system
$(f(x),D)$ and $B_{r_{1}}(x_{1})$ and $B_{r_{1}}(x_{2})$ are
respectively the attractor $A(x_{1})$ and the Siegel's disk
$SI(x_{2})$.}\\

{\bf Case of $p>2$, $|2a|_p=|b|_p$ and $|2a-b|_p=|a|_p$}.

In this case we have $|2a|_p=|a|_p=|b|_p=|2a-b|_p$.
By similar argument as in the previous case we can show that
$B_{r_{0}}(x_{1})=A(x_{1})$.

We now consider $x_{2}$ and its Siegel disk. The condition for $r$
in this case is $\Gamma 2$ and  it suggests, that
$$\max_{n\geq 2}\left(|a-b|_pr\right)^{n-1}<\left(\frac{|a|_p}{|a-b|_p}\right)^{2}.$$
From this we deduce, that if $r>\frac{1}{|a-b|_p}$, then
$\left(|a-b|_pr\right)^{n-1}$ is unbounded and, thus - we must
necessarily require that $r\leq\frac{1}{|a-b|_p}$ for the maximum
to have a finite magnitude. Then, if $r=\frac{1}{|a-b|_p}$ holds,
$\left(|a-b|_pr\right)^{n-1}$ is essentially equal to $1$ for
$\forall n\geq 2$ which can only be possible when the right-hand
side of the conditional inequality
$\left(\frac{|b|_p}{|a-b|_p}\right)^{2}>1$ $\Leftrightarrow$
$|a-b|_p<|b|_p$. If, instead we have that $|a-b|_p=|b|_p$ then
should be $r<\frac{1}{|a-b|_p}$ and then the maximal value is
attained for $n=2$ and is equal to $|b|_pr$. In this case $r$ must
satisfy the inequality $r<\frac{|b|_p}{|a-b|_p^{2}}$. Now, this is
consistent with the assumption $r<\frac{1}{|a-b|_p}$ since
$\frac{|b|_p}{|a-b|_p^{2}}=\frac{|b|_p}{|a-b|_p}\frac{1}{|a-b|_p}<\frac{1}{|a-b|_p}$.
This uncertainty gives rise to the following two possibilities

\textbf{\tt Possibility 1.} $|a-b|_p=|b|_p$ in this case the
inequality $r<\frac{1}{|a-b|_p}$ should hold for
$B_{r}(x_{2})\subseteq SI(x_{2})$ to be true.

\textbf{\tt Possibility 2.} $|a-b|_p<|b|_p$ in this case the
inequality $r\leq\frac{1}{|a-b|_p}$ should hold for
$B_{2}(x_{2})\subseteq SI(x_{2})$ to be true.

Note, that in both cases $B_{\frac{1}{|a-b|_p}}(x_{2})\subseteq
SI(x_{2})$ but $\bar{B}_{\frac{1}{|a-b|_p}}(x_{2})\nsubseteq
SI(x_{2})$ since the point
$x_{0}=-\frac{1}{a}\in\bar{B}_{\frac{1}{|a-b|_p}}(x_{2})$ but
$f(x_{0})=x_{2}$ and, thus, $x_{0}\notin SI(x_{2})$. This shows,
that $B_{\frac{1}{|a-b|_p}}(x_{2})=SI(x_{2})$ in both of the
possible case above.

Thus we have proved the following

\textbf{Theorem 3.7.} \textit{When $p>2$, $|2a|_p=|b|_p$ and
$|2a-b|_p=|a|_p$, then the points $x_{1}$ and $x_{2}$ are
respectively the attracting point and the indifferent point of the
dynamical system $(f(x),D)$ and $B_{\frac{1}{|a|_p}}(x_{1})$ and
$B_{\frac{1}{|a-b|_p}}(x_{2})$ are respectively the attractor
$A(x_{1})$ and the Siegel's disk $SI(x_{2})$.}

By similar argument as above the same theorem can be proved for the case of $p=2$ and
$|b|_p=|a|_p$.

\subsection{$x_{2}$ is an attracting point}

Here, as was demonstrated, we have three distinct cases, namely:

\begin{enumerate}
    \item The case of $p=2$ and $|2a|_p>|b|_p$.
    \item The case of $|2a|_p=|b|_p$ and $|2a-b|_p<|a|_p$.
    \item The case of $p=2$ and $|2a|_p<|b|_p<|a|_p$.
\end{enumerate}

{\bf Case of $p=2$ and
$|2a|>|b|_p$.}

In this case we have  $|2a-b|_p=|2a|_p=\frac{1}{2}|a|_p$, $|a|_p>|b|_p$ and
$|a-b|_p=|a|_p$.

By similar reasoning as in the previous subsection we can prove that
$B_{r_{1}}(x_{1})=A(x_{1})$.

Now, consider $x_{2}$ and its attractor $A(x_{2})$. The condition
for $r$ in this case is $\Gamma 3$ - it suggests, that
$$\max_{n\geq 2}\frac{|a|_p}{|b|_p}\left(|b|_pr\right)^{n-1}<1.$$
Whence we find that $B_{r_{1}}(x_{2})\subseteq A(x_{2})$.

Now, if $|x-x_{2}|_p\geq r_{0}$ then $|x|_p\geq \dsf{1}{|b|_p}$.
According to the equality $\Delta 3$ have
$$
|f(x)-x_{2}|_p\geq\dsf{|a|_p}{|b|_p}|x-x_{2}|_p\geq |x-x_2|_p.
$$
Thus $x$ does not belong to the $A(x_{2})$.

If, $r_{1}<r<r_{0}$, i.e. $|x-x_{2}|_p=r$, this means that
$|x|_p=r>|x_{2}|_p$ and $\left|x+\frac{1}{a}\right|_p=|x|$
then, modifying equality $\Delta 3$ we derive the following
$$
|f(x)-x_{2}|_p|P|_p=\dsf{|a|_p}{|b|_p}|x-x_{2}|_p^{2}
$$
and hence
$$
|f(x)-x_{2}|_p=|a|_p|x-x_{2}|_p^{2}.
$$
Now, $|f(x)-x_{2}|_p>|a|_p|x-x_{2}|_p\frac{1}{|a|_p}=|x-x_{2}|_p$.
Thus, such $x$'s do not lie in the $A(x_{2})$.

Now suppose that $|x-x_{2}|_p=r_{1}=\frac{1}{|a|_p}$ then this implies
$$
|ax-bx-1|_p=1.\eqno(3.4)
$$
We have to consider a case $|ax|_p=1$, since $|ax|_p<1$ means $x\in A(x_1)$.

On the other hand from $\Delta 3$  we get
$$
|f(x)-x_{2}|_p=\bigg|x+\dsf{1}{a}\bigg|_p.\eqno(3.5)
$$
Let us assume that
$$
|f(x)-x_{2}|_p<|x-x_{2}|_p=\dsf{1}{|a|_p}.
$$
it then follows from (3.5) that $|ax+1|_p<1$, whence $|ax-1|_p<1$ since
$p=2$. From (3.4) we find $|bx|_p=1$, which contradicts to our assumption
$|a|_p>|b|_p$. Hence
$|f(x)-x_{2}|_p=\dsf{1}{|a|_p}$. So there are two possible cases
(i) $|f(x)|_p<\dsf{1}{|a|_p}$ and (ii) $|f(x)|_p=\dsf{1}{|a|_p}$.
If (i) is satisfied then $x\in A(x_1)$, if (ii) is satisfied then
by similar argument as above we have $|f^2(x)-x_{2}|_p=\dsf{1}{|a|_p}$.
Consequently, we obtain in this setting $x\notin A(x_2)$. Thus
$A(x_2)=B_{r_1}(x_2)$.

Hence we have proved the following

\textbf{Theorem 3.8.} \textit{When $p=2$ and $|2a|_p>|b|_p$ the
points $x_{1}$ and $x_{2}$ are the attracting points of the
dynamical system $(f(x),D)$ and $B_{r_{1}}(x_{1})$ and
$B_{r_{1}}(x_{2})$  are
the attractor respectively  $A(x_{1})$ and $A(x_{2})$.}\\

By similar argument as above the same theorem can be proved for the case of
$|2a|_p=|b|_p$ and $|2a-b|_p<|a|_p$.\\

{\bf Case of $p>2$ and $|2a|_p=|b|_p$ and
$|2a-b|_p<|a|_p$}.

In this case we have $|a-b|_p=|a|_p$ and $|a|_p=|b|_p$, so
$r_1=r_0$.

For the fixed point $x_{1}$ and the $A(x_{1})$ the $\Gamma 1$
condition requires that
$$\max_{n\geq 2}(|b|_pr)^{n-1}<1
$$
which implies $r<\dsf{1}{|b|_p}=r_0$. Hence, for all $r<r_0$
we have $B_{r}(x_1)\subseteq A(x_{1})$.
By the similar argument as in the previous subsections
we can prove that  $B_{r}(x_1)=A(x_{1})$.

We turn our attention to the consideration of the point $x_{2}$
and its attractor $A(x_{2})$. The condition for $r$ in this case
is $\Gamma 3$ and it implies that
$B_{r_{0}}(x_{2})\subseteq A(x_{2})$.

Now, if $|x-x_{2}|_p>r_{0}$ then we have $|x-P|_p=|x|_p$. By means of
$\Delta 3$ we get
$|f(x)-x_2|_p=|x-x_2|_p$ hence $x\notin A(x_2)$, here we have used
$|x+\frac{1}{a}|_p=|x|_p$.

Let us turn to a case $|x-x_2|_p=r_0$. In this case we have to consider
$|x|_p=r_0$ since $|x|_p<r_0$ implies $x\in A(x_1)$. It is clear that
$P\in S_{r_0}(x_2)$, therefore $|x-P|_p\leq r_0$. Assume
$|x-P|_p=r', \ r'<r_0$. Then from $\Delta 3$ we infer
$$
|f(x)-x_2|_p=\bigg|x+\frac{1}{a}\bigg|_p\dsf{r_0}{r'}.
\eqno(3.6)
$$
Now we show that r.h.s. of (3.6) cannot be less than $r_0$. Suppose contrary,
which means
$$
\bigg|x+\frac{1}{a}\bigg|_p<{r'}.
\eqno(3.7)
$$
Using strong triangle property and
$$
\bigg|\dsf{1}{a}+P\bigg|_p=\bigg|\dsf{a-b}{ab}\bigg|_p=r_0
$$
we have
$$
|x-P|_p=\bigg|\left(x+\dsf{1}{a}\right)-\left(P+\dsf{1}{a}\right)\bigg|_p=r_0.
$$
The last equality contradicts to $|x-P|<r_0$. Hence in the considered
case $x\notin A(x_2)$.

Now assume $|x-P|_p=r_0$. Then $\Delta 3$ implies
$$
|f(x)-x_2|_p=\bigg|x+\frac{1}{a}\bigg|_p.
$$
Let us suppose that $\bigg|x+\frac{1}{a}\bigg|_p<{r_0}$. Then using the inequality
$$
\bigg|\dsf{1}{a}+x_2\bigg|_p=\bigg|\dsf{2a-b}{a(a-b)}\bigg|_p<r_0
$$
and property 3.1 we obtain
$$
|x-x_2|_p=\bigg|\left(x+\dsf{1}{a}\right)-\left(x_2+\dsf{1}{a}\right)\bigg|_p<r_0.
$$
which contradicts to $|x-x_2|_p=r_0$, hence $x\notin A(x_2)$. So we
get $A(x_2)=B_{r_0}(x_2)$.

\textbf{Theorem 3.9.} \textit{When $p>2$ and $|2a|_p=|b|_p$ and
$|2a-b|_p<|a|_p$ the points $x_{1}$ and $x_{2}$ are the attracting
points of the dynamical system $(f(x),D)$ and $B_{r_{0}}(x_{1})$
and $B_{r_{0}}(x_{2})$  are
the attractor respectively  $A(x_{1})$ and $A(x_{2})$.}\\

By similar argument as above the same theorem can be proved for the case of
$p=2$ and $|2a|_p<|b|_p<|a|_p$.

\section{$p$-adic dynamical system of the form
$f(x)=\frac{x^{2}}{bx+a}$ is not ergodic}

In this section we will show that the considered dynamical system
is not ergodic with respect to the Haar measure.

Let us first consider a dynamical system of the form
$$
f(x)=\dsf{x^2}{bx+1}, \ \ \ b\in\Q_p, \ b\neq 0 \eqno(4.1)
$$ in $\Q_p$.  It is easy to see that $x_2=\dsf{1}{1-b}$ is a fixed point for (4.1).
A question about ergodicity of the considered system arises when
the fixed point $x_2=\dsf{1}{1-b}$ is indifferent. This leads us
to the condition $|b|_p<1$ (see subsection 3.2). In what follows
we will suppose that $p>2$. We consider our system on the sphere
$S_{\rho}(x_2)$ with $\rho<1$. The indifference of the point $x_2$
implies that $f(S_{\rho}(x_2))=S_{\rho}(x_2)$.

From Chapter 27 of \cite{E} we infer the following

{\bf Lemma 4.1.}{\it  For every ball $B_{r}(a)\subset S_{\rho}(x_2)$, $r<\rho$
the following equality holds
$$
f(B_{r}(a))=B_{r}(f(a))
$$
}

Consider a measurable space $(S_\rho(x_2),{\cal B})$, here ${\cal
B}$ is the algebra of generated by clopen subsets of
$S_\rho(x_2)$. Every element of ${\cal B}$ is a union of some
balls $B_{r}(a)$, $r=p^{-l}<\rho$. A measure $\mu:{\cal B}\to \R$
is said to be {\it Haar measure} if it is defined by
$$
\mu(B_{p^{-l}}(a))=\frac{1}{q^l},
$$
here $q$ is a prime number.

From Lemma 4.1 we conclude that $f$ preserves the measure $\mu$,
i.e.
$$
\mu(f(B_{r}(a)))=\mu(B_{r}(a))\eqno(4.2).$$

Now we are going to prove that $f$ is not ergodic on
$S_\rho(x_2)$, i.e. we will find a subset $A$ of $S_\rho(x_2)$
such that $f(A)=A$ and $0<\m(A)<1$ ( we refer a reader for some
preliminaries on ergodic theory to \cite{W}).

{\bf Lemma 4.2.}{\it  For every ball $B_{r}(y)\subset
S_{\rho}(x_2)$, $r<\rho$ the following equality holds
$$
B_r(y)\cap B_r(f(y))=\emptyset.
$$
}

{\bf Proof.} It is enough to prove the equality $|f(y)-y|_p=\rho$.
Note that $|y-P|_p=|P|_p=\dsf{1}{|b|_p}$, here $P=-\dsf{1}{b}$.
Now using $\Delta 2$ we have
$$
|f(a)-y|_p|y-P|_p=\dsf{1}{|b|_p}|y|_p|y-x_2|_p,
$$
which implies the required equality. The proof is completed.

{\bf Lemma 4.3.} {\it The following equality holds for every $x\neq P$
$$
(f^2(x)-x)(x-P)=\dsf{1-b}{b}\left[\dsf{f(x)}{(x-P)b}\left(x+1\right)+x\right](x-x_2),
$$
here as before $f^2(x)=f(f(x))$. }

The proof is based on the equalities $\Delta 2$ and $\Delta 3$ (see Proposition 3.2).

{\bf Lemma 4.4.} {\it For $r_0=\rho|b|_p$ and $y\in S_{\rho}(x_2)$
the following equality holds
$$
B_{r_0}(y)=B_{r_0}(f^2(y)).
$$
}

{\bf Proof.} In order to prove the required equality it is enough
to show that $|f^2(y)-y|_p\leq r_0$. It is clear that
$B_{r_0}(y)\subset S_{\rho}(x_2)$. As $B_{|1-b|}(x_2)$ is the
Siegel disk (see Theorem 3.7) and $S_\rho(x_2)\subset
B_{|1-b|}(x_2)$, we have
$$
|f(y)-x_2|_p=|y-x_2|_p=\rho
$$
whence we infer $|f(y)|_p=1$, since $|x_2|_p=1$. Using Lemma 4.3
we conclude
$$
|f^2(y)-y|_p=\dsf{1}{|y-P|_p}\left|\dsf{f(y)}{(y-P)b}\left(y+1\right)+y\right|_p|y-x_2|_p.
\eqno(4.3)
$$
Keeping in mind (4.2) we get
$$
\left|\dsf{f(y)}{(y-P)b}\left(y+1\right)+y\right|_p\leq 1,
$$
which implies with (4.3) that
$$
|f^2(y)-y|_p\leq |b|_p\rho.
$$

This completes the proof.

The proved Lemma 4.4 and Lemma 4.2 mean that the following set
$$
A=B_{r_0}(y)\cup B_{r_0}(f(y))
$$
is invariant with respect to $f(x)$ and $0<\mu(A)\neq \mu(S_{\rho}(x_2))$.
Hence, $f(x)$ is not ergodic on $S_{\rho}(x_2)$ with respect to the Haar measure.

Thus we have the following

{\bf Theorem 4.5.} {\it The dynamical system (4.1) is not ergodic
on $S_\rho(x_2)$ with respect to the Haar measure.}

{\bf Remark 4.1.}  In \cite{GKL} it has been proved that a
dynamical system $g(x)=x^2$ is ergodic on $S_{\rho}(1)$ with
respect to Haar measure when $p=3$. Our dynamical system (4.1)
when $b=0$ becomes $g(x)$. But Theorem 4.5 shows that if we
slightly perturb the $g(x)$ the obtained dynamical system is not
ergodic for any $p>2$.

{\bf Remark 4.2.} Let consider a dynamical system of the form
$$
h(x)=\dsf{x^2}{bx+a},  \ \ \ 0\neq |b|_p<1,\ a\in Q_p,
\eqno(4.4)
$$
on $S_{\rho|a|}(\tilde x_2)$, here as before $\rho<1$ and $\tilde x_2=|a|x_2$ - fixed point
of (4.4). From the results of subsection it is easy to check that $S_{\rho|a|}(\tilde x_2)$
is a Siegel disk of (4.4). On $S_{\rho|a|}(\tilde x_2)$ we consider Haar measure.
Now define a function
$$
S(x)=ax, \ \ \ \ x\in\Q_p.
$$
It is evident that $SfS^{-1}(x)=h(x)$ for all $x\in\Q_p$. Hence, according to
Theorem 4.5 we infer that dynamical system (4.4) is not ergodic.\\

{\bf Acknowledgements}  One of authors (F.M.) thanks CNR for
providing financial support and II Universita di Roma "Tor
Vergata" for all facilities, where the final part of the paper was
completed. The work is also partially supported by Grants
$\Phi$-1.1.2, $\Phi$-2.1.56  of CST of the Republic of Uzbekistan.
%\newpage

\edc